\documentclass[fleqn,11pt]{wlscirep}
\usepackage[utf8]{inputenc}
\usepackage[T1]{fontenc}
\usepackage{bm}
\usepackage{float}
\usepackage{comment}
\usepackage{cleveref}
\usepackage[dvipsnames]{xcolor}

\newcommand{\abs}[1]{\left|#1\right|}

\usepackage{ulem}
\usepackage{xcolor,tikz}
\newcommand\delete{%
 \bgroup\markoverwith
  {\textcolor{orange}{\pgfsetfillopacity{0}\rule[-0.5ex]{2pt}{10pt}\pgfsetfillopacity{1}}%
   \textcolor{orange}{\llap{\rule[0.5ex]{2pt}{1pt}}}%
  }%
  \ULon}
\usepackage{hyperref}

\title{Discovery of Unstable Singularities}

\author[1, 2]{Yongji Wang}
\author[3]{Mehdi Bennani}
\author[3]{James Martens}
\author[3]{S\'{e}bastien Racani\`{e}re}
\author[3]{Sam Blackwell}
\author[3]{Alex Matthews}
\author[3]{Stanislav Nikolov}
\author[1, 4]{Gonzalo Cao-Labora}
\author[5]{Daniel S. Park}
\author[3, $^{\dagger}$]{Martin Arjovsky}
\author[3, $^{\dagger}$]{Daniel Worrall}
\author[3, $^{\dagger}$]{Chongli Qin}
\author[3, $^{\dagger}$]{Ferran Alet}
\author[3, $^{\dagger}$]{Borislav Kozlovskii}
\author[3, $^{\dagger}$]{Nenad Toma\v{s}ev}
\author[3]{Alex Davies}
\author[3]{Pushmeet Kohli}
\author[1,*]{Tristan Buckmaster}
\author[3,*]{Bogdan Georgiev}
\author[6,*]{Javier G\'{o}mez-Serrano}
\author[3,*]{Ray Jiang}
\author[2,*]{Ching-Yao Lai}

\affil[1]{New York University, Department of Mathematics, New York, NY 10012, USA}
\affil[2]{Stanford University, Department of Geophysics, Stanford, CA 94305, USA}
\affil[3]{Google DeepMind, London, N1C 4DJ, UK}
\affil[4]{\'Ecole Polytechnique F\'ed\'erale de Lausanne, Institute of Mathematics, Lausanne, VA 1015, Switzerland}
\affil[5]{Google DeepMind, New York, NY 10011, USA}
\affil[6]{Brown University, Department of Mathematics, Providence, RI 02912, USA}

\affil[$^{\dagger}$]{These authors have comparable contributions}
\affil[*]{Corresponding author: buckmaster@cims.nyu.edu}
\affil[*]{Corresponding author: bogeorgiev@google.com}
\affil[*]{Corresponding author: javier\_gomez\_serrano@brown.edu}
\affil[*]{Corresponding author: rayjiang@google.com}
\affil[*]{Corresponding author: cyaolai@stanford.edu}
\affil[*]{All corresponding authors contributed equally and are listed in alphabetical order}

\begin{abstract}
Whether singularities can form in fluids remains a foundational unanswered question in mathematics. This phenomenon occurs when solutions to governing equations, such as the 3D Euler equations, develop infinite gradients from smooth initial conditions. Historically, numerical approaches have primarily identified stable singularities. However, these are not expected to exist for key open problems, such as the boundary-free Euler and Navier-Stokes cases, where unstable singularities are hypothesized to play a crucial role.  Here, we present the first systematic discovery of new families of unstable singularities. A stable singularity is a robust outcome, forming even if the initial state is slightly perturbed. In contrast, unstable singularities are exceptionally elusive; they require initial conditions tuned with infinite precision, being in a state of instability whereby infinitesimal perturbations immediately divert the solution from its blow-up trajectory. In particular, we present multiple new, unstable self-similar solutions for the incompressible porous media equation and the 3D Euler equation with boundary, revealing a simple empirical asymptotic formula relating the blow-up rate to the order of instability.  Our approach combines curated machine learning architectures and training schemes with a high-precision Gauss-Newton optimizer, achieving accuracies that significantly surpass previous work across all discovered solutions. Crucially, for specific solutions (the Córdoba-Córdoba-Fontelos (CCF) stable and 1st unstable solutions), we reach near double-float machine precision, attaining a level of accuracy constrained only by the inherent round-off errors of the GPU hardware. This level of precision meets the stringent requirements for rigorous mathematical validation via computer-assisted proofs. This work provides a new playbook for exploring the complex landscape of nonlinear partial differential equations (PDEs) and tackling long-standing challenges in mathematical physics.
\end{abstract}
\begin{document}

\flushbottom
\maketitle

\thispagestyle{empty}

Understanding whether the partial differential equations governing fluid motion can develop singularities from smooth initial conditions -- where mathematical solutions exhibit infinite gradients or velocities in finite time -- represents one of the most enduring challenges in mathematics. This fundamental question spans a broad class of nonlinear fluid equations, including the incompressible 3D Euler equations \cite{Euler:euler-equations} (first formulated in 1752), the Boussinesq equations \cite{Majda-Bertozzi:vorticity-incompressible-flow} which model atmospheric flows, and the incompressible porous media (IPM) equation \cite{Bear:dynamics-porous-media}, describing fluid motion through porous materials like soil or rock. Crucial insights into these dynamics are also gained through the study of simplified analogues, such as the 1D Córdoba-Córdoba-Fontelos (CCF) model \cite{Cordoba-Cordoba-Fontelos:CCF-model}.

Among these systems, the 3D Navier-Stokes equations, the viscous analogue of the 3D Euler equations, represent perhaps the most famous case. The question of finite-time blow-up for Navier-Stokes is recognized as one of the seven Millennium Prize Problems \cite{Fefferman:clay-statement} and the Smale Problems \cite{Smale:problem-list} (the $15^{th}$), with profound implications for mathematics and physics. We refer to the surveys \cite{Hou:survey-blowup-euler-ns,Gibbon:3d-euler-survey} for developments on this problem. If singularities can develop from smooth initial conditions, the equations predict physically impossible outcomes, such as infinite velocity gradients, in finite time. This would signal a fundamental breakdown in the predictive power of the equations. 

A crucial aspect of singularity formation is stability. A singularity is considered stable if it reliably emerges from a range of initial conditions. Conversely, an unstable singularity requires infinitely precise initial conditions; any slight perturbation causes the system to deviate from the blow-up trajectory. It is hypothesized that singularities in foundational, boundary-free problems like the 3D Euler and Navier-Stokes equations must be unstable. While the presence of boundaries can stabilize solutions, developing the techniques to discover and resolve unstable singularities, even in systems {\it with} boundaries, is a critical intermediate step toward tackling the boundary-free challenge. Furthermore, it is expected that highly unstable solutions are better candidates for persisting when transitioning from idealized equations to more realistic, viscous ones, making their discovery essential for resolving these fundamental questions.

Efforts to resolve the problem of singularity formation have spanned both numerical simulations and classical analysis. In terms of numerical candidates to Euler blow-up, several works have considered the possibility of colliding dipoles\cite{PumirSiggia,kerr,Hou2022}. By performing space-time simulations up to times very close to the singularity, Luo and Hou\cite{Luo-Hou:singularities-euler-3d} presented compelling numerical evidence for finite-time blow-up for the axisymmetric 3D Euler equations if one introduces a cylindrical boundary.

An emerging approach combines high-precision numerics with rigorous mathematical analysis in the form of a computer-assisted proof (CAP) via the use of interval arithmetic\cite{GomezSerrano:survey-cap-in-pde}. The self-similar profile corresponding to the Luo-Hou scenario was first found numerically using physics-informed neural networks (PINNs)\cite{Raissi-Perdikaris-Karniadakis:pinn-first-paper,Karniadakis-Kevrekidis-Lu-Perdikaris-Wang-Yang:pinn-nature-reviews} by Wang et al. \cite{Wang-Lai-GomezSerrano-Buckmaster:pinn-selfsimilar-boussinesq} and proven rigorously using a CAP by Chen and Hou\cite{Chen-Hou:stable-nearly-self-similar-blowup-boussinesq-part-i,Chen-Hou:stable-nearly-self-similar-blowup-boussinesq-part-ii} (c.f. Wang et al.\ \cite{2025arXiv250619243W}). A key ingredient in this hybrid approach\cite{Costin-Tanveer:analytical-approximation-blasius,Elgindi-Pasqualotto:invertibility-boussinesq-cap,Castro-Cordoba-GomezSerrano:global-smooth-solutions-sqg} is the use of self-similar coordinates. These coordinates rescale space and time around the point of blow-up \cite{Buckmaster-CaoLabora-GomezSerrano:implosion-compressible,Eggers-Fontelos:self-similarity,Elgindi:singularities-C1a-euler-R3}. This removes the difficulty of tracking rapidly diverging physical quantities by transforming the problem from a time-evolving simulation into finding a stationary, smooth self-similar profile.

This work addresses the following critical, previously unsolved challenge. Previous frameworks utilizing numerical tools such as B-splines on non-uniform meshes \cite{Luo-Hou:singularities-euler-3d}, finite elements \cite{Guillod-Sverak:numerics-nonuniqueness-ns} or spectral methods \cite{Lushnikov-Silantyev-Siegel:self-similar-gclm} excel at resolving stable blow-up solutions, yet they may face significant limitations in capturing unstable ones. Performing time-stepping simulations (even self-similar) to identify unstable solutions requires extreme precision to overcome the exponential instabilities arising from being on the unstable manifold. This limitation is more notable when those solutions are quantized (as opposed to a continuum of solutions in terms of a continuous parameter) and continuation methods cannot be applied.
Indeed, in the case of colliding dipoles, instabilities present preclude the use of such numerical techniques as a basis for a mathematical construction of blow-up. In order to render the problem accessible by such numerics, Hou \cite{Hou:potential-singularity-euler} added viscous terms to stabilize the blow-up; however, due to the singular nature of such terms, it is not clear how the solutions found relate to the original equations considered. The critical challenge, therefore, has been developing a computational framework capable of discovering and accurately resolving these elusive, unstable solutions to the precision required for mathematical proof.

\begin{figure}[t]
	\centering
    \includegraphics[
        width=\linewidth,
        trim={0mm 147mm 0mm 0mm},
        clip
    ]{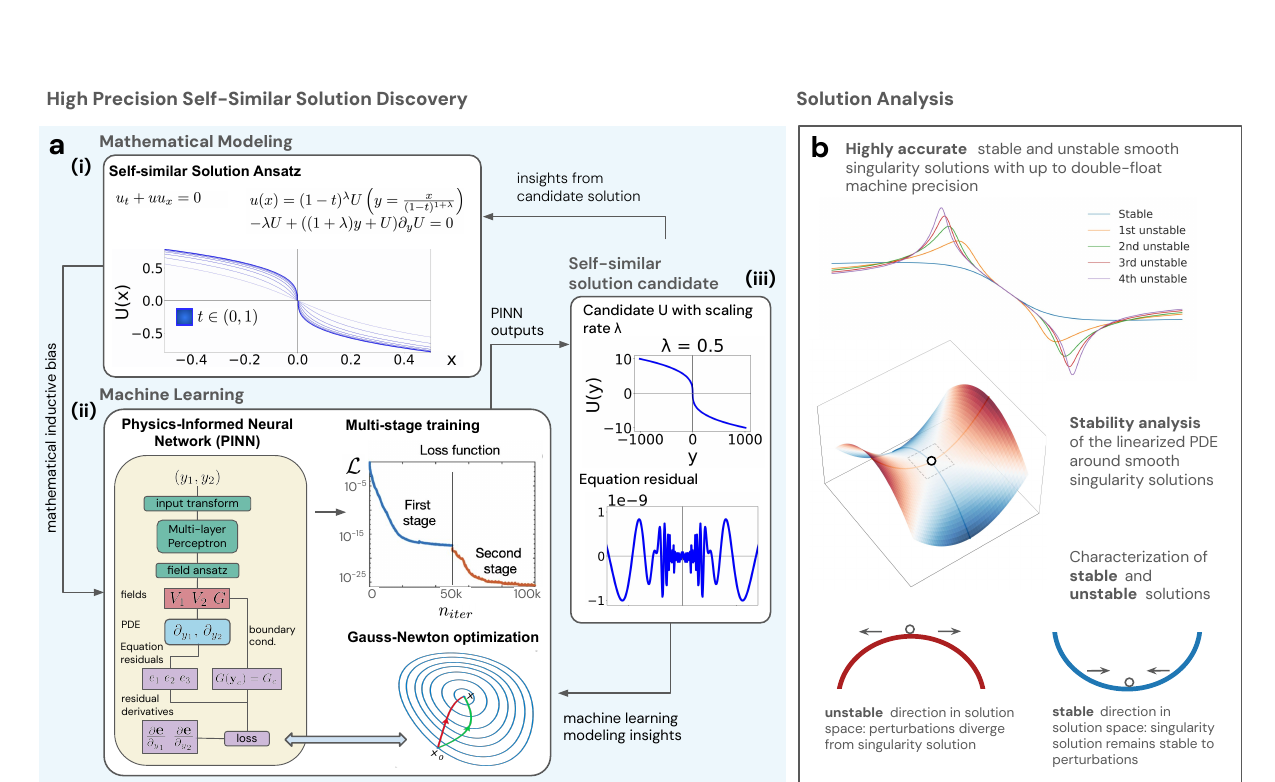}
	\caption{\small \textbf{Research flowchart.} Our research methodology consists of two main stages. \textbf{a) Solution discovery.} We start with a candidate solution that allows us to search for self-similar spatial profiles of blow-up solutions characterized by a self-similar scaling rate $\lambda$, illustrated in \textbf{(i)} with Burgers' Equation. We then refine our machine learning pipeline \textbf{(ii)} and solution accuracy using an iterative approach. Empirical results from candidate solutions \textbf{(iii)} and their accuracy guide the mathematical modeling and neural network architecture. Mathematical modeling in turn guides the inductive biases built into the network architecture, such as input coordinate transforms and the shaping of the output fields. We use a Physics-informed Neural Network (PINN) with a Gauss-Newton optimiser and a multi-stage refinement training scheme to generate a highly accurate candidate solution while finding the right scaling rate $\lambda$. \textbf{b) Solution analysis.} We analyze the stability of a high-precision candidate solution by linearising the PDE around it. We discover unstable modes, along which any slight perturbation causes the system to deviate from the blow-up trajectory. Thus, we characterize the degree of stability and confirm the discovery of highly accurate stable and unstable singularities.
    }
	\label{fig:roadmap}
\end{figure}

\begin{figure}[t]
	\centering
	\includegraphics[width=0.98\linewidth]{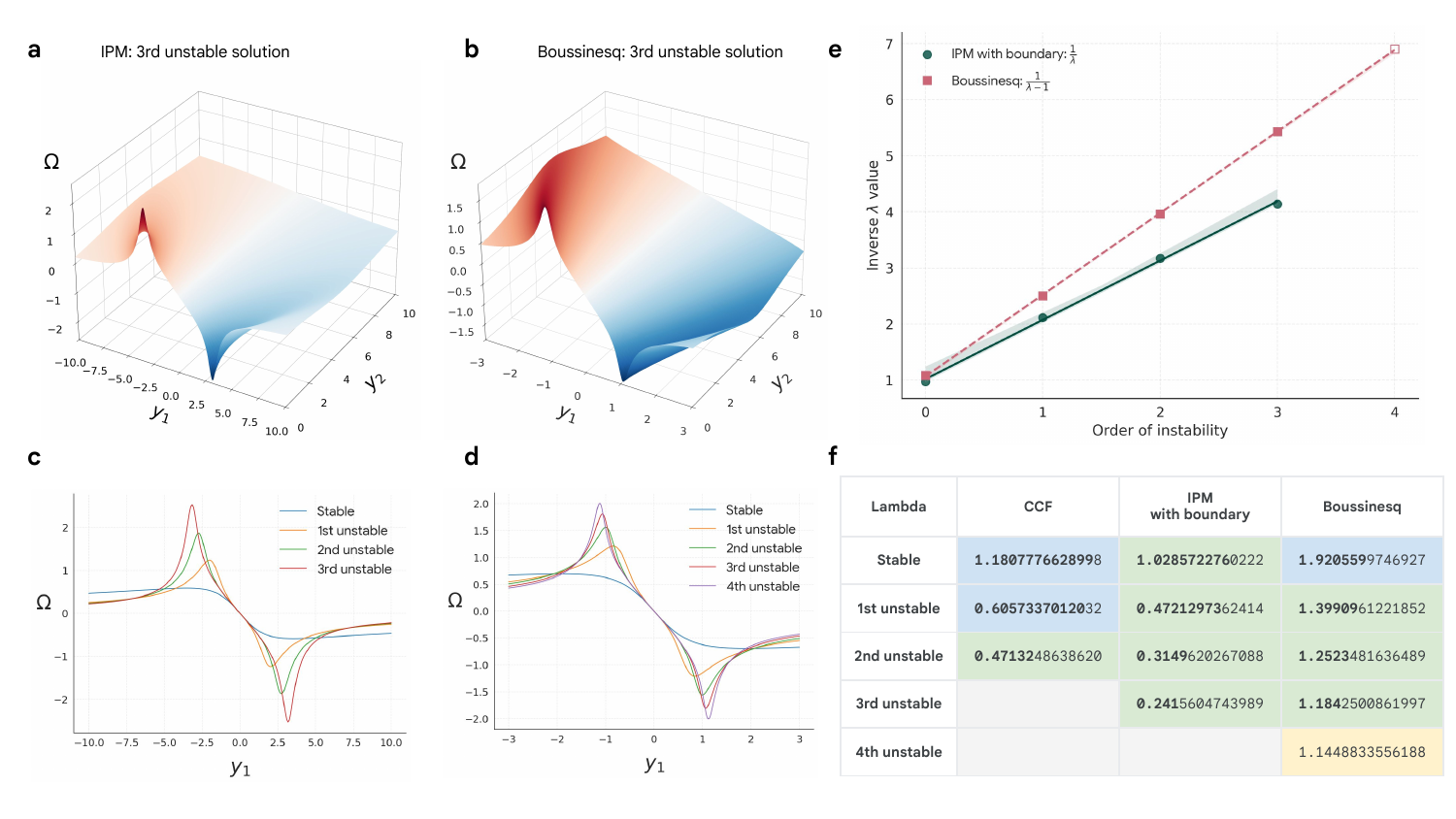}
	\caption{\small \textbf{Self-similar singularities to IPM and Boussinesq.}  \textbf{a, b,} Spatial profile of the vorticity $\Omega$ for the third unstable solutions of the IPM (in \textbf{a}) and Boussinesq equations (in \textbf{b}) near the origin. \textbf{c, d,} The cross section of the vorticity profiles along the $y_1$ axis, namely $\Omega(y_1, 0)$ for the stable to third unstable solutions of the IPM (in \textbf{c}) and Boussinesq equations (in \textbf{d}). \textbf{e,} Inverse scaling rate versus the instability order, the plot reveals a linear trend for both systems (shaded area denotes one standard deviation). \textbf{f,} The scaling parameters of the unstable solutions discovered in the paper (green and yellow) and the previously discovered solutions (blue). The significant digits of the scaling parameters are marked in bold-face. The 4th unstable solution of the Boussinesq equation is un-validated.}
	\label{fig:blowup}
\end{figure}

In this work, we introduce such a framework (see Figure \ref{fig:roadmap}), combining curated neural network architectures and training schemes with a high-precision Gauss-Newton optimizer. This approach enables the first systematic discovery of new families of unstable singularities for the IPM and Boussinesq equations, revealing a simple empirical asymptotic formula relating the blow-up rate to the order of instability (see Figure \ref{fig:blowup}). We achieve a level of accuracy that significantly surpasses previous work across all discovered solutions; crucially, for specific CCF solutions, we reach near double-float machine precision, a level constrained only by the inherent round-off errors of the GPU hardware. 

\section*{Discovery of Unstable Singularities}

In this study, we discover unstable singularities in three canonical fluid systems: the CCF equations\cite{Cordoba-Cordoba-Fontelos:CCF-model}, the 2D IPM equations, and the 2D Boussinesq equations. In the presence of a boundary, the 2D Boussinesq equations are mathematically analogous (up to exponentially small terms) to the incompressible 3D Euler equations under the assumption of axial symmetry.

To discover singularities, each equation is reformulated in self-similar coordinates -- a time-dependent rescaling of space parameterized by a scaling parameter $\lambda$ that governs the blow-up spatiotemporal structure (see  Methods \hyperref[methods:self-similar]{Self-similar Coordinates}). By construction, finding a self-similar singularity entails the task of identifying $\lambda$ values admitting smooth solutions in those coordinates. We say such values of $\lambda$ are \textit{admissible}. The self-similar parameter $\lambda$ dictates the nature of the singularity (see Methods \hyperref[methods:lambda]{Identifying admissible $\lambda$}). A smaller value of $\lambda$ corresponds to a more unstable solution and, consequently, a more severe singularity with a faster blow-up rate, which allows it to overcome greater dissipative effects (e.g.\ going from Euler to Navier-Stokes).

We quantify the accuracy of the discovered numerical solutions using the maximum residual, which is defined as the maximum absolute value of the equation residual over the entire spatial domains when a candidate solution is substituted into the governing equations. To establish a meaningful and comparable metric across different solutions, it is essential to fix the scale (normalization). The maximum (absolute) residual, calculated across all equations and a dense spatial validation grid, therefore provides a standardized measure of how accurately the solution satisfies the system. Beyond evaluating the residual, we employ other validation methods. Crucially, we analyze the spectrum of the linearisation of the governing equations around the candidate solution, checking for expected properties related to the stability of the solution, such as the number of unstable modes. Once a numerical solution and its corresponding scaling factor $\lambda$ is determined, we validate the significant digits of $\lambda$ by finding the minimum amount of perturbation that can be applied to $\lambda$ for which the solution quality is unchanged. We discuss these validation methods in detail in Methods and the Supplementary Information.\footnote{The Supplementary Information to the paper will be uploaded in a forthcoming version of the manuscript.}

The CCF equations can be viewed as a simplified model for the 3D Euler equations and the Navier-Stokes equation. An important open question\cite{Silvestre-Vicol:transport-equation-nonlocal-drift}\cite[Open Problem 1]{Kiselev:regularity-blowup-active-scalars} regarding the CCF equations is whether singularities still occur in the presence of fractional dissipation $(-\Delta)^{\alpha/2}$ for $\frac12 \leq \alpha \leq 1$. This problem involves the fundamental balance between nonlinear amplification (which causes blow-up) and dissipative mechanisms (which regularize the solution). The first unstable blow-up profile for CCF was discovered by Eggers and Fontelos\cite{eggers_2020}, and recently confirmed using PINN-based methods by Wang et al. \cite{Wang-Lai-GomezSerrano-Buckmaster:pinn-selfsimilar-boussinesq} with $\lambda_1 \approx 0.6057$, suggesting that singularities exist for any $\alpha \leq 1/(1+\lambda_1) \approx 0.623$, ultimately culminating in a mathematically rigorous proof (manuscript in preparation). The challenge of finding higher unstable solutions for the CCF equations lies in the extremely sharp gradients of their profiles, often believed to trap the optimization in local minima. In this work, we improve the stable and first unstable singularities of Wang et al.\cite{Wang-Lai-GomezSerrano-Buckmaster:pinn-selfsimilar-boussinesq} to achieve PDE residuals of order $O(10^{-13})$ for both solutions. We also discover the second unstable profile with $\lambda_2 = 0.4703$ with residuals of order $O(10^{-7})$. This discovery increases the critical $\alpha$ for singularity formation to $\alpha \leq 1/(1+\lambda_2) \approx 0.68$, advancing our understanding of dissipative thresholds in fluid systems and shedding light on an old and fundamental open question. 

The IPM equations describe density evolution driven by incompressible flows governed by Darcy’s law under gravity \cite{Bear:dynamics-porous-media}. While prior literature provides no evidence of self-similar solutions in IPM\cite{cordobagancedooriveipm2007,cordobazoroaipm},
we identify the stable and  \textit{three} unstable singularities with residuals whose orders of magnitude range from $O(10^{-11})$ to $O(10^{-8})$.

For the Boussinesq equations, we are able to confirm the stable solution and discover \textit{three} distinct unstable singularities with residuals whose orders of magnitude range from $O(10^{-8})$ to $O(10^{-7})$. Furthermore, we have identified a candidate for the fourth unstable solution. This solution, however, is not resolved to the same level of accuracy as the others; in particular, we do not yet have a reliable evaluation of the accuracy of its corresponding $\lambda$ value. We refer to the Supplementary Information for more details regarding this candidate. This marks the first discovery of any smooth unstable self-similar solution to any unforced incompressible fluid equation, leading to a finite time blow-up. We refer to Methods and Supplementary Information for the presentation and analysis of all solutions.

For every unstable solution discovered for the CCF, IPM and Boussinesq equations, we conduct a stability analysis by linearizing the governing equations around the solution. For the $n$-th unstable solution discovered, we find $n$ unstable modes that respect the same symmetry assumptions as the solution. Since such unstable modes specify directions in which the solution can be pushed to become more stable, this suggests that the family of solutions discovered so far are complete within the range of admissible $\lambda$ values considered.

One proposed approach to addressing the Clay Millennium Prize Problem involves finding a sequence of self-similar profiles, corresponding to a sequence of self-similar scaling parameters $\lambda_n$, for the 3D Euler equations, each with an increasing number of unstable directions. It is expected that higher-order instability makes it more feasible to treat viscosity as a perturbative error, as the timescales in which the viscous effects come into play are much slower than the ones leading to the finite time singularity. This phenomenon has been observed in the compressible Navier-Stokes equations and the defocusing Nonlinear Schr\"odinger (NLS) equation\cite{Merle-Raphael-Rodnianski-Szeftel:implosion-ii,Merle-Raphael-Rodnianski-Szeftel:implosion-nls}.

Based on the discovery of families of unstable singularities we have made, we are able to extrapolate an empirical rule for admissible $\lambda$ values of the Boussinesq/Euler equations and the IPM equations. For the Boussinesq/Euler equations, we find that $\lambda_n \sim 1/(1.4187n + 1.0863) + 1$, where $\lambda_n$ indicates the $\lambda$ value for the $n$-th unstable singularity. For IPM with boundary, the empirical relation of $\lambda$ for higher-order unstable modes is $\lambda_n \sim 1/(1.1459n + 0.9723)$. Both asymptotic expansions fit well for larger values of $n$. For the CCF equations, a clear asymptotic relationship has not yet been identified from the solutions discovered so far, making a reliable extrapolation premature. These formulae provide candidate values for admissible $\lambda$ for higher-order unstable modes, which can guide future studies of singularity hierarchies in these equations. Moreover, they can make numerical discovery of new solutions easier by providing high-quality initialization values for $\lambda$, making thus the distance between the initial value at training and the final value smaller, yielding faster and more reliable convergence.

\section*{Targeted PDE solution discovery with PINNs}

In this study, all unstable singularities are discovered by modeling the solutions with physics-informed neural networks (PINNs), which has emerged as a new tool for solving differential equations \cite{Raissi-Perdikaris-Karniadakis:pinn-first-paper,Karniadakis-Kevrekidis-Lu-Perdikaris-Wang-Yang:pinn-nature-reviews}. Wang et al. \cite{Wang-Lai-GomezSerrano-Buckmaster:pinn-selfsimilar-boussinesq} have demonstrated the potential of using PINNs for discovering stable singularities in fluids. PINNs have also been applied to explore unstable numerical solutions, for instance, in 2D cavity flow problems~\cite{wang2023solution, zou2025learning}. So far, the accuracies commonly achieved by PINNs (e.g., an error of $O(10^{-3})$ for Boussinesq in Wang et al. \cite{Wang-Lai-GomezSerrano-Buckmaster:pinn-selfsimilar-boussinesq}) have been insufficient for the numerical solutions to be used for a rigorous proof.

The application of a PINN to a problem is characterized by three core components: the neural network architecture, the loss function, and the optimization algorithm. Substantial research has been focused on improving the performance of PINNs by advancing these components. This includes developments in network design (such as activation function selection \cite{jagtap2020adaptive, wang2023learning} and configuration \cite{jagtap2021extended, moseley2023finite}), sophisticated strategies for weighting the loss function \cite{mcclenny2020self, wang2021understanding, wang2022and, wang2022respecting}, the adoption of more effective optimization methods \cite{muller2023achieving,Jnini-Vella-Zeinhofer:gauss-newton-natural-gradient-pinn} and the employment of hierarchical training schemes \cite{aldirany2024multi, Wang-Lai:multistage-nn}. 

Given that there is recent concern about overoptimism in using PINNs as general-purpose PDE solvers \cite{mcgreivy2024weak}, we stress that PINNs are not used for this purpose in this work. Rather, the goal is to discover solutions to differential equations that have not been found before. Here, PINNs are utilized to parameterize the particular solution of partial differential equations we wish to discover, and the pipelines for training the models are designed and constructed with the intention of making the network satisfy all the mathematical assumptions and constraints imposed by the problem at hand.

\begin{figure}[t]
	\centering
	\includegraphics[width=0.95\linewidth]{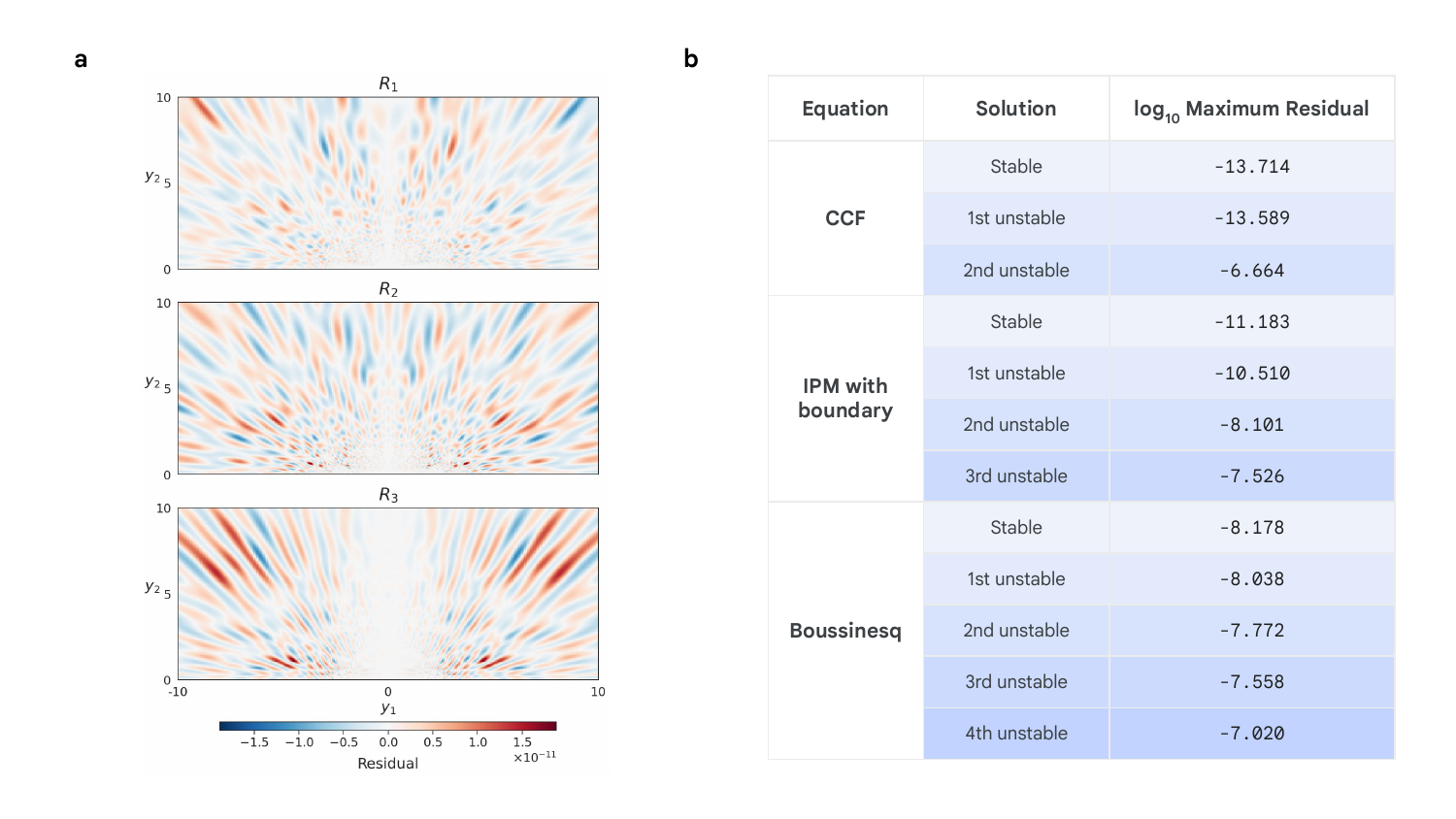}
	\caption{\small \textbf{Solution accuracy.} \textbf{a,} Spatial profiles of the equation residual for the 1st unstable solution to the IPM equation after multistage training.  \textbf{b,} List of the log-10 maximum equation residuals for the best solutions discovered for the CCF, IPM and Boussinesq equations.}
	\label{fig:all_residuals}
\end{figure}

Since our goal is discovery, we do not claim improvement in efficiency or speed against any benchmarks. The validation metrics we present only concern the quality of the solutions that have been discovered. In figure \ref{fig:all_residuals}\textbf{b}, we present the best maximum residuals we have been able to achieve for all solutions. The maximum residuals generally range from $O(10^{-8})$ to $O(10^{-7})$, achieving multiple orders of magnitude improvement in this metric compared to the literature \cite{Wang-Lai-GomezSerrano-Buckmaster:pinn-selfsimilar-boussinesq,2025arXiv250619243W}. For the CCF equations, we are able to approach precision limited only by the double-floating point round-off errors ($O(10^{-13})$), a level of accuracy currently required for enabling subsequent computer-assisted proofs.

\subsection*{Incorporating mathematical structure into a neural network}

We represent solutions as smooth functions parameterized by neural networks. This allows us to directly embed known mathematical properties of the target solution into the network architecture. This embedding acts as a strong inductive bias, guiding the optimization toward mathematically relevant solutions and away from trivial or degenerate ones.

We enforce constraints derived from the governing equations through architectural design. Symmetries and periodicity are imposed via input transformations. Infinite domains, inherent in self-similar problems, are handled using coordinate transformations that compactify the space\cite{Lushnikov-Silantyev-Siegel:self-similar-gclm} (see Methods \hyperref[methods:solution_modeling]{Solution Modeling}). Furthermore, we utilize tailored "solution envelopes" -- multiplicative factors applied to the network output -- to enforce behaviors such as asymptotic decay at infinity and local series expansions near the origin. Crucially, smoothness requirements at the origin can also be used to derive an analytical relationship for the self-similar parameter $\lambda$ (see \cite{Dahne-GomezSerrano:highest-wave-burgershilbert} for an analogous application). This provides a robust alternative to learning $\lambda$ as a standard trainable parameter \cite{Wang-Lai-GomezSerrano-Buckmaster:pinn-selfsimilar-boussinesq}.

While formal analysis provides essential guidance, it often relies on assumptions (e.g. selecting one branch from several formal possibilities) that are not guaranteed a priori. We overcome this through a feedback loop between numerical experiments and mathematical analysis (see Figure \ref{fig:roadmap}). Initial experiments reveal hidden structures which are then incorporated back into the network architecture a posteriori. For example, experiments for the IPM and Boussinesq equations suggest that the solutions vanish faster near the origin than required by minimal symmetry assumptions. Incorporating this insight allows us to explicitly factor this vanishing behavior out from the solution representation. We then are able to reformulate the governing equations for the remaining component, which significantly improves optimization stability by structuring the terms the network must satisfy (see Methods and Supplementary Information).

By systematically embedding these mathematical features, both those known a priori and those discovered iteratively, we enhance the optimization landscape. This integration of mathematical structure accelerates convergence and ultimately enables the discovery of previously inaccessible high-order unstable solutions.

\begin{figure}[t]
	\centering
	\includegraphics[width=0.95\linewidth]{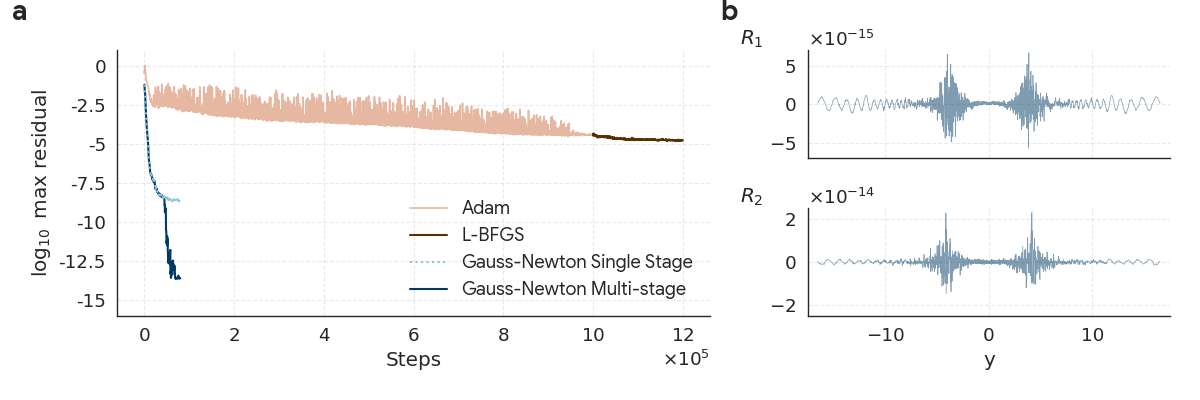}
	\caption{\small \textbf{High-precision training.} \textbf{a,} The log-10 maxium residual value of the solution evaluated on the evaluation grid for the CCF first unstable solution when training with different optimizers. The Gaussian-Newton optimizer shows better performance and significantly faster convergence than the classical Adam + L-BFGS optimizers. Multi-stage training is able to further improve the results by five orders of magnitude. \textbf{b,} The two residuals of the CCF first unstable solution obtained by multi-stage training plotted against the self-similar spatial coordinate $y$.}
	\label{fig:ablation_test}
\end{figure}

\subsection*{High-precision training}

Improving the training procedure is a crucial element to achieving the extreme precision required for discovering unstable singularities. We find the two key components to attaining high-precision training to be the utilization of a full-matrix Gauss-Newton optimizer, and the employment of multi-stage training.

As presented in figure \ref{fig:ablation_test}, we find standard gradient optimization methods, e.g., Adam or L-BFGS, to be insufficient for producing high quality solutions to the self-similar equations, c.f. McGreivy and Hakim\cite{mcgreivy2024weak}. Recent works have explored the use of stronger second-order optimizers, which involve approximations to the curvature matrix to remain computationally feasible for large networks, for related tasks\cite{muller2023achieving,Jnini-Vella-Zeinhofer:gauss-newton-natural-gradient-pinn,dangel2024kronecker}.

In contrast, we adopt a {full-matrix} Gauss-Newton (GN) approach\cite{nocedal1999numerical}, to optimize our neural networks\cite{martens2020new}. Since no structural approximations to the GN matrix are made, we are able to compute an essentially unbiased stochastic estimate of it. While generally computationally infeasible for most applications of neural networks, this approach is viable for us due to the relatively small size of our networks, whose number of parameters range from thousands to at most tens of thousands. 

We efficiently estimate the GN matrix using an inexpensive, unbiased rank-1 estimator, which we average across iterations using an exponential moving average to reduce variance. To ensure robust convergence without manual tuning, we integrate this approach with automated methods\cite{martens2015optimizing} for determining nearly optimal per-iteration learning rates and momentum coefficients. This (nearly) fully automated scheme, implemented using the open-source library \texttt{kfac-jax}\cite{kfac-jax2022github}, allows us to rapidly iterate on modeling choices without expensive hyperparameter sweeps. In \Cref{fig:ablation_test} we demonstrate the efficacy of this optimizer, which is able to discover solutions with maximum residuals of $O(10^{-8})$ within 50k iterations ($\approx$ 3 A100 GPU hours), compared to standard optimization techniques.

Using a better optimizer alone does not get us to the precision of solutions we would like to ultimately achieve. To that end, we employ multi-stage training\cite{Wang-Lai:multistage-nn}, in which a second neural network is optimized to efficiently correct the remaining high frequency error of the first network that has already been trained. The output of the first and second networks are combined to produce the solution. We utilize multi-stage training for the stable and first unstable solutions of the CCF and IPM equations. For the solutions of the CCF equations, we are able to improve the maximum residuals by five orders of magnitude down to $O(10^{-13})$, a level sufficient for rigorous mathematical validation based on a CAP. See \cref{fig:ablation_test} for the validation curve for multi-stage training applied to the first unstable solution of CCF, as well as the spatial distribution of the residuals of the discovered solution. For the IPM stable and first unstable solutions, we are able to achieve maximum residuals of $O(10^{-11})$ and $O(10^{-10})$ respectively through multi-stage training (Fig. \ref{fig:all_residuals}\textbf{a}).

\section*{Implications for Future Studies}

Since their introduction\cite{Raissi-Perdikaris-Karniadakis:pinn-first-paper}, PINNs have been widely studied in science and engineering \cite{Karniadakis-Kevrekidis-Lu-Perdikaris-Wang-Yang:pinn-nature-reviews}. Here, we demonstrate their advantage as a specialized discovery tool for finding new and highly accurate numerical solutions that advance mathematical understanding. Rather than proposing PINNs as general-purpose PDE solvers, we leverage their designability to seek specific, elusive solutions -- such as the unstable singularities detailed in this paper -- that pose considerable challenges for discovery. Our approach is able to reduce PDE residuals by multiple orders of magnitude compared to those reported in the literature, achieving double-float machine precision for solutions of the CCF equation. This level of accuracy proves to be sufficient for utilizing the numerical solution for computer-assisted proofs for this specific problem\footnote{Manuscript in preparation.}, although the required precision for such proofs is generally problem-dependent. We hope our results would encourage new numerical designs and discovery methods to develop for similar applications across the computational community. Our results highlight the importance of combining domain-specific mathematical insights with modern computational tools to tackle problems that have resisted purely analytical or purely numerical methods. This interdisciplinary approach opens new possibilities for exploring complex solution landscapes in nonlinear PDEs. However, significant challenges remain. A key challenge is to discover self-similar solutions to the incompressible 3D Euler equations (or similar equations such as IPM) in the absence of boundary. To address this problem, we will need to further understand the qualitative properties of such solutions that can guide the design of PINN architecture and training pipelines.

\bibliography{sample,references}

\newpage
\section*{Methods}
\label{methods}

\phantomsection
We aim to discover self-similar blow-up solutions to a system of fluid equations with high numerical accuracy using Physics-Informed Neural Networks (PINNs). The governing equations that define the problem describe the time evolution of physical quantities of the fluid such as fluid density or velocity. Such fields can be written as mappings of the form $\phi_A(x,t): \mathbb{R}^d \times \mathbb{R} \rightarrow \mathbb{R}$ for fields labeled by $A = 1, 2, \cdots, N_F$, where $N_F$ is the number of field components of the problem. We use $d$ to denote the spatial dimension involved. We can omit the field index and collect the fields into a vector:
\begin{equation}
    \phi(x,t): \mathbb{R}^d \times \mathbb{R} \rightarrow \mathbb{R}^{N_F} \,.
\end{equation}
We use $N_E$ to denote the number of governing equations.

\subsection*{Self-similar Solutions}
\label{methods:self-similar}
We seek solutions whose initial profile evolves in a self-similar manner and blows up at time $t=1$. To do so, we introduce the self-similar coordinate $y$ and assume that the behavior of the fields with respect to time is fixed to have a power law behavior. In other words, we assume that there exists a self-similar scaling parameter $\lambda$, and for each field with label $A$, a constant exponent $k_A(\lambda)$ such that:
\begin{equation}
\label{eq:self-similar}
    \phi_A(x,t) = (1-t)^{k_A(\lambda)} \Phi_A(y) \quad \text{with}
    \quad
    y=(1-t)^{-(1+\lambda)} x 
\end{equation}
solves the governing equations. When $d > 1$, it is understood that both $y$ and $x$ are vectors related by the time-dependent multiplicative factor. The exponents $k_A(\lambda)$ are linear functions of $\lambda$ assigned in a way such that a power of $(1-t)$ can be factored out of each governing equation. Put differently, $k_A$ is determined so that the governing equations, labeled by index $k = 1, 2, \cdots, N_E$, factor into the form $(1 -t)^{f_k(\lambda)} \cdot \mathcal{E}_k (\Phi(y), \lambda) = 0$ upon plugging in the ansatz of \cref{eq:self-similar}. The original dynamical equations are now converted into the steady equations
\begin{equation}
\label{eq:self-similar-static}
 \mathcal{E}_k (\Phi(y), \lambda) = 0\,,
\end{equation}
satisfied by the smooth self-similar fields
\begin{equation}
    \Phi(y) : \mathbb{R}^d \rightarrow \mathbb{R}^{N_F} \,,
\end{equation}
which only depend on the spatial dimension.

The self-similar scaling parameter $\lambda$ plays a crucial role in characterizing the self-similar solution, as the governing equations themselves have an explicit dependency on the parameter $\lambda$. In fact, for many examples, there is evidence, and in few cases proof\cite{eggers2008}, that smooth self-similar solutions only exist for isolated values of $\lambda$.  As noted in the main text, we say such values of $\lambda$ are \textit{admissible}.

Through the transformation of \cref{eq:self-similar} we have converted the original problem to finding an admissible value of $\lambda$ and the corresponding smooth self-similar profile $\Phi(y)$ satisfying the static equations \cref{eq:self-similar-static}. We model the self-similar profile with neural networks and apply specialized optimization methods to solve the self-similar equations.

\subsection*{Solution Modeling and Equation Factorisation}
\label{methods:solution_modeling}
A challenging aspect of solving the self-similar blow-up equations is that they are defined over infinite or semi-infinite spatial domains. Numerically solving these equations in a truncated finite domain introduces significant (truncation) error. In order to overcome this issue, we design a coordinate transformation that maps the entire infinite domain to a finite one\cite{grosch1977numerical,lushnikov2017new,boyd1987spectral}. For ease of exposition, we restrict our attention to two-dimensional problems in our presentation, although analogous devices are employed for the one-dimensional case. More details specific to each problem can be found in the Supplementary Information.

A key mathematical constraint we need to satisfy is the asymptotic behavior of the profile at infinity. The asymptotic behavior is dominated by some of the linear terms in the equation, e.g.\ for the Boussinesq equation, we expect the vorticity to behave like $\Omega = \text{curl}({\bf U}) \sim  \abs{y}^{-\alpha}$ for the self-similar spatial variable $y=(y_1,y_2)$ and exponent $\alpha = 1/(1+\lambda)$.

Based on the power-law decay, we introduce the new coordinate system $z = (q, \beta)$ whose components are
\begin{equation}
    q = (\sqrt{1+\abs{y}^2})^{-\alpha} = (1+\abs{y}^2)^{-1/(2(1+\lambda))} \,, \quad
    \beta = y_2/\sqrt{1+y_1^2 + y_2^2} \,.
\end{equation}
With these new coordinates, we translate a problem defined on an unbounded domain to that of a bounded domain, since the map $y = (y_1, y_2) \mapsto z = (q, \beta)$ maps the $y$-halfplane to a subset of $(0, 1] \times (-1, 1)$. We pass the pair $(q, \beta)$ as the input to the neural networks. This allows the networks to learn the profile of the solution within a normalized, finite domain, where its functional variation is expected to be simpler than in the original coordinates.

This coordinate transformation introduces an implicit constraint. Since $q(y_1, y_2)$ and $\beta(y_1, y_2)$ are both even functions of the coordinate $y_1$, any standard neural network receiving them as input will produce an output that is also inherently even in $y_1$. This is not a problem if the field we intend to model has this symmetry, but is an issue if the behavior is expected to be different. We thus construct a structured solution ansatz for each physical quantity based on their symmetry properties. Rather than having the network learn the solution directly, we have it learn a simpler, even function, which is then multiplied by a prefactor designed to enforce the correct symmetry and asymptotic behavior.

Again, we use the vorticity, $\Omega$, from the Boussinesq equation as an example. For this problem, we look for an $\Omega$ which is an odd function of $y_1$ and decays at infinity with the power-law exponent $\alpha$. We formulate the expression for $\Omega$ as:
\begin{equation} \label{eq:ansatz_vorticity}
    \Omega(y_1, y_2) = \underbrace{ \left( \frac{y_1}{\sqrt{1+y_1^2+y_2^2}} \right) }_{\text{Odd Symmetry}} \cdot  \mathrm{NN}_{\Omega}[q, \beta]  \cdot \underbrace{ q(y_1, y_2) }_{\text{Asymptotic Decay}}
\end{equation}
where the prefactor is a simple odd function that imposes the correct symmetry on the entire expression. The neural network output, $\mathrm{NN}_{\Omega}[q, \beta]$, learns the remaining even part of the spatial profile of the solution in the transformed coordinates. The final multiplication by \textbf{$q(y_1, y_2)$} explicitly enforces the required power-law decay rate at infinity.

This construction is generalized for each variable in the problems we solve, allowing us to hard-code both symmetry and asymptotic boundary conditions directly into the model architecture. This strategy offloads the burden of learning global physical constraints from the optimizer and shifts it onto the network design itself. As a result, the network can focus its capacity on learning the localized features of the solution.

We use multi-layer perceptrons (MLPs) with $\tanh$ activations as the neural network component. We choose to use an MLP in this context because we want a versatile, hierarchical, nonlinear function approximator, building on the success of Wang et al. \cite{Wang-Lai-GomezSerrano-Buckmaster:pinn-selfsimilar-boussinesq} in this specific context.  We use small neural networks with thousands to tens of thousands of parameters. The output of the fields must capture fluctuations of multiple orders of magnitude. Thus the final layer of the network is set to be an exponential or exponential-adjacent activation layer.

Meanwhile, we find that the governing equations of the self similar solutions of \cref{eq:self-similar-static} can factorize in a useful way. More precisely, writing the neural network component of the field as
\begin{equation}
    \hat{\Phi}_\theta: \mathbb{R}^d \rightarrow \mathbb{R}^{N_F} \qquad
    \text{where} \quad
    \hat{\Phi}_{A,\theta} = \text{NN}_{A,\theta} (z(y)) \,,
\end{equation}
we find that the governing equations factor in the form:
\begin{equation}
    \mathcal{E}_k (\Phi(y), \lambda) = F_k(y) \cdot \mathcal{R}_k (\hat{\Phi}_\theta(y), \lambda) \,,
\end{equation}
such that $F_k$ is well defined in the domain of the fields. We add the $\theta$ subscript here to make the dependence of $\hat{\Phi}$ on the neural network parameters $\theta$ explicit. $F_k$ may be the identity for some equations, in which case $\mathcal{E}_k$ and $\mathcal{R}_k$ would be identical.

We factor the effect of $F_k$ out from the governing equations and denote the remaining factors $\mathcal{R}_k (\hat{\Phi}_\theta(y), \lambda)$, "residuals" of the equations. By doing so, we dispose of the components of the governing equations that are mechanically determined by the assumptions we make about the solution or by the structure of the equation. Now the machine learning problem becomes the problem of finding parameters of the neural networks $\hat{\Phi}_\theta$ that minimize the size of the residuals $\mathcal{R}_k (\hat{\Phi}_\theta(y), \lambda)$. The precise definition of the equation residuals for the CCF, IPM and Boussinesq equations can be found in the Supplementary Information.

\subsection*{Training Loss}
The goal of training is to make the magnitude of residuals $\mathcal{R}_k$ as small as possible across the entire domain of the fields. Thus, the key component of the loss is the zeroth-derivative, or "d0," loss:
\begin{equation}
    \mathcal{L}_\text{d0} (Y_0) = {1 \over N_\text{d0} N_E} \sum_{(y, \rho) \in Y_0} \rho \sum_{k=1}^{N_E} \left[ \mathcal{R}_k(\hat{\Phi}_\theta(y), \lambda) \right]^2\,.
\end{equation}
The dependence of the loss on the sampling batch $Y_0$ is made explicit in this expression. An element of a batch is given by the pair $(y, \rho)$, where $y$ is the sampling point, and $\rho$ is the loss weight assigned to that point.

Only training with the d0 loss can result in a spiky solution, that only has small residuals in a small neighborhood around the d0 training points. To make the residuals as close to a constant as possible, we also utilize higher derivative losses:
\begin{align}
    \mathcal{L}_\text{d1}(Y_1) &= {1 \over N_\text{d1} N_E d} \sum_{(y, \rho) \in Y_1} \rho \sum_{k=1}^{N_E} \sum_{i=1}^d \left[ \partial_{y_i} \mathcal{R}_k(\hat{\Phi}_\theta(y), \lambda) \right]^2\,, \\
    \mathcal{L}_\text{d2}(Y_2) &= {1 \over N_\text{d2} N_E d^2} \sum_{(y, \rho) \in Y_2} \rho \sum_{k=1}^{N_E} \sum_{i, j=1}^{d} \left[ \partial_{y_i} \partial_{y_j} \mathcal{R}_k(\hat{\Phi}_\theta(y), \lambda) \right]^2\,,    
\end{align}
where we have used $i$ and $j$ to index the spatial dimensions. The "equation loss" is expressed as a linear combination of these losses:
\begin{equation}
    \mathcal{L}_\text{equation} = 
    c_\text{d0} \cdot \mathcal{L}_\text{d0} +
    c_\text{d1} \cdot \mathcal{L}_\text{d1} +
    c_\text{d2} \cdot \mathcal{L}_\text{d2} \,.
\end{equation}
We note that each loss has its own assigned training batch.

Without additional constraints imposed, the neural network is bound to produce a constant zero function. Furthermore, there are rescaling symmetries that can act on a solution to produce another solution. In order to avoid such trivial solutions and fix such symmetries, we add normalization conditions on the network where we impose constraints on the field values or its derivatives at certain points. To do so, a point $y_\text{norm}$ is chosen at a high gradient region at which the value of a specific field $\Phi_{A_0}$ is required to take some value $C_\text{norm}$. Also, the same field is required to have vanishing value at a set of points $Y_\infty$ far away from the origin. A typical "data loss" incorporating these basic constraints takes the form:
\begin{align}
    \mathcal{L}_\text{data} = 
    \hat{c}_\text{norm} \left(\Phi_{A_0}(y_\text{norm}) - C_\text{norm}\right)^2 
    + \hat{c}_\infty \sum_{y \in Y_\infty} \Phi_{A_0}(y)^2 \,.
\end{align}
The constraint points are fixed throughout training. For certain problems, we may introduce constraints for more than one field or derivatives of the field. We also may enforce bounds by acting on a constraint with a relu function before squaring it to add to the data loss. Details on the exact data loss used for each problem can be found in the Supplementary Information. The training loss is obtained by adding the data loss to the equation loss:
\begin{equation}
    \mathcal{L}_\text{train} = \mathcal{L}_\text{equation} + \mathcal{L}_\text{data} \,.
\end{equation}

\subsection*{Validation Metrics for Solutions}

In order to judge the quality of individual solutions during and after training, we compute the maximum d$n$ residuals of the solutions on a dense validation grid. While the grid may be constructed to span an arbitrary large range within the domain, for example, using the compactified coordinates, practically, the points with large residuals are focused around small regions with high gradients not far from the origin. Thus, for the purpose of evaluation, using a grid that densely samples those regions turns out to be more important. The grids used for evaluating solutions of each problem are specified in detail in the Supplementary Information.

Given a validation grid, the maximum d$n$ residual is defined to be the maximum absolute value of the $n$-derivative of the residuals, taken along the grid points, the possible derivatives, and the governing equations. More formally,
\begin{equation}
    \text{max d0 residual}_\text{grid} (\hat{\Phi}_\theta, \lambda) =
    \underset{\substack{k \in [N_F] \\y \in\text{grid}}}{\text{max}} \left|
    \left[ \mathcal{R}_k (\hat{\Phi}_\theta(y), \lambda)\right] \right|\,,
\end{equation}
while more generally for $n \geq 1$,
\begin{equation}
    \text{max d$n$ residual}_\text{grid} (\hat{\Phi}_\theta, \lambda) =
    \underset{\substack{k \in [N_F]\\ i_1, \cdots, i_n \in [d]\\y \in\text{grid}}}{\text{max}}
    \left| \left[ \left( \prod_{m=1}^n \partial_{i_m} \right) \mathcal{R}_k (\hat{\Phi}_\theta(y), \lambda)\right] \right| \,.
\end{equation}

\subsection*{Training Collocation Points}

An important element of the training process is the selection of collocation points. Since the residuals of the equations can be computed for any point in space, we must specify how collocation points in space are to be sampled to produce a training batch for the model. We have two challenges. The first is to sample points that have high information content. For example, the solution behaves trivially when we move far away from the origin, while the cancellation of terms in the governing equation are highly non-trivial in regions where the field values and their derivatives are large. Meanwhile, we would like to avoid overfitting, where the model ends up only focusing on lowering the residuals in select regions and choosing to neglect others.

To address both of these challenges, we employ two different sampling methods. The first is location based sampling, where points are sampled from a fixed region in space with an empirical measure. How the regions are selected, how the sampling measure is set, and what loss weights are assigned to the samples are determined empirically and explained in detail in the Supplementary Information.

The next is adaptive sampling, where points are globally sampled but in a way so that points where the model currently does not perform well, i.e., where the residuals of the governing equations are large, are sampled more frequently. In practice, a global grid of points are defined over the compactified coordinates $z$. The residuals $\mathcal{R}_k$ and their derivatives are evaluated over the grid points, and used to define the adaptive sampling weight. For example, for the Boussinesq equation, we define the sampling weight to be proportional to the fourth power of the squared sum of the residuals. After the grid points are sampled, a small perturbation may be applied to the grid point so that a neighboring point is sampled.

For each loss $\mathcal{L}_\text{d0}$, $\mathcal{L}_\text{d1}$ and $\mathcal{L}_\text{d2}$, we define a sampling scheme for the training batches $Y_0$, $Y_1$ and $Y_2$ passed to the loss. The sampling scheme consists of defining point groups and how many samples should taken from each of the point groups. For each group, we designate a sampling method. If location based sampling is being used, we specify the region of space to be sampled from, what the sampling measure is and what the loss weight assignment should be. Typically, a global region, as well as regions close to the origin or peaks of the solutions are utilized. If adaptive sampling is being used, we specify the sampling weight as a function of the residuals and their derivatives, as well as the loss weights for the points. The sampling schemes, regions, measures and loss weights used for the CCF, IPM and Boussinesq equations are presented in full detail in the Supplementary Information.

Diminishing the residuals on a set of collocation points turns out to be a difficult task. Thus resampling the training points very frequently does not lead to meaningful training. We therefore train with a given set of collocation points for an extended number of steps, and only resample them after an order of thousands of steps have elapsed.

\subsection*{Second-order Optimizer}
The PINN training objective function $\mathcal{L}(\theta, y)$ can be written as $h(\theta, y) = \sum_i c_i \|f_i(\theta, y)\|^2$ for nonlinear vector-valued functions $f_i$, scalar constants $c_i > 0$, parameter vector $\theta \in \mathbb{R}^n$, and a vector $y$ representing all of the collocation points. Given an objective of this form, we can apply the classical Gauss-Newton optimization method\cite{nocedal1999numerical}, which is a preconditioned gradient descent of the form $\theta_{t+1} = \theta_t - (G(\theta_t, y) + \gamma_t I)^{-1} \nabla h (\theta, y)$, where $G$ is the Gauss-Newton matrix and $\gamma$ is a ``damping" coefficient.

Unlike other recent works\cite{muller2023achieving,Jnini-Vella-Zeinhofer:gauss-newton-natural-gradient-pinn,dangel2024kronecker,rathore2024challenges,wang2025gradient} applying second-order optimization to PINNs, which use approximations that simplify the structure of $G$, we compute an essentially unbiased stochastic estimate of the full $G$ matrix. In particular, we estimate $G$ using an exponential moving average (across optimizer iterations) of an unbiased rank-1 estimator of $G$ that can be computed from a single matrix-vector product with the Jacobian of the $f_i$'s. This full-matrix approach is feasible due to the small size of our networks. 

Because of estimation variance and slight bias of the exponential moving average, we cannot simply use a constant learning rate of $1$ (which is viable when $G$ is computed precisely), and instead adopt the scheme described by Martens and Grosse\cite{martens2015optimizing} to compute "optimal" learning rate and momentum values in closed form at each iteration. This approach increases the per-iteration cost of the method by factor of about 2 to 3 (excluding the cost of matrix inversion, which is non-dominant), but frees us from having to tune learning rate schedules via trial and error and expensive sweeps, and in particular gives us a schedule customized to our training problem.

The optimization algorithm involves the damping coefficient $\gamma$, which controls how much we trust the local quadratic approximation of the objective underlying our second-order optimizer. How $\gamma$ is adjusted throughout optimization is crucial to the robust performance of second-order methods. Depending on the problem, we either use the classical Levenberg-Marquardt method \cite{levenberg1944method,marquardt1963algorithm} or employ custom schedules for adjusting the value of $\gamma$.

This optimizer can be fully realized using a particular configuration of the open source \texttt{kfac-jax} library\cite{kfac-jax2022github}, which is a general framework supporting different 2nd-order optimizers including, but not restricted to, K-FAC. For a fully detailed description of our optimization method and its implementation, see the Supplementary Information.

\subsection*{Identifying Admissible $\lambda$ Values}
\label{methods:lambda}
The scaling parameter $\lambda$ plays a central role in studying self-similar solutions. We operate with the assumption that there exist discrete admissible values of $\lambda$ that admit smooth solutions to the self-similar equations, which we present empirical evidence for later on. Thus identifying admissible values of $\lambda$ is a key element to discovering solutions to the self-similar equations.

In this work, we simultaneously optimize $\lambda$ as we train the parameters of the neural network. Three different methods for $\lambda$ optimization are employed. The combination of $\lambda$ optimization methods used for each problem and solution can be found in the Supplementary Information.\\
\\
{\bf 1. Gradient training:} The simplest way to optimize $\lambda$ is to introduce it as a trainable parameter and optimize it along with the other neural network parameters\cite{Wang-Lai-GomezSerrano-Buckmaster:pinn-selfsimilar-boussinesq}. While gradient training is able to successfully identify admissible values of $\lambda$ for stable blow-up solutions in the CCF, IPM, and Boussinesq problems, we find it to be unreliable for discovering unstable, high-order solutions.\\
\\
{\bf 2. $\lambda$ Inference:} We use the $\lambda$ inference method (c.f.\ Prat-Colomer \cite{PratColomer:master-thesis} in the context of ML and rotating solutions of the Euler equation), where we infer the value of $\lambda$ best suited for describing the current numerical solution from the governing equations. Assuming the smoothness of the solution, we take the first derivative of the nonlinear self-similar equation for the given problem and evaluate it at the origin. Requiring that this value vanishes allows us to derive an explicit analytical relationship between $\lambda$ and field values at the origin $(y_1, y_2) = (0, 0)$. For example, the relation derived for the Boussinesq problem is that $\lambda = -3-2\partial_{y_1}U_1(0, 0)$, where $\partial_{y_1}U_1(0,0)$ is the velocity gradient at the origin. We then implement this relationship as a direct analytical update rule for $\lambda$ during training, bypassing the need for an optimization-based search.\\
\\
{\bf 3. Funnel Inference:} Another method we introduce for optimizing $\lambda$ is funnel inference. This method is based on the principle that the governing equation can be perfectly satisfied only at admissible values of $\lambda$. Thus in the neighborhood of an admissible value $\lambda^*$ of the self-similar scaling parameter, the smallest achievable maximum residual as a function of $\lambda$ should have a funnel-like shape, the bottom of which is positioned at $\lambda^*$. If we consider this residual with sign, in other words, define the residual with maximum absolute value to be $r(\lambda)$, $r(\lambda)$ should be approximately linear near $\lambda^*$ where it has zero value. The idea of funnel inference is to use the secant method to find the zero of this function $r(\lambda)$.

Since $r(\lambda)$ is not directly measurable, we find a proxy for $r(\lambda)$. To do so, for a given constant $\lambda$ value, we proceed to train the neural network with this fixed $\lambda$ value until the solution has relaxed sufficiently. Operating with the assumption that this signal can be accurately captured by examining a small neighborhood of the origin, we measure the residuals of the equation and identify the value of the residual with maximum absolute value in this region to be $\hat{r}(\lambda)$. Then $\lambda$ is updated to a new value, and the training process repeats. The update rule is given as follows. Given the history of $\lambda$ values $\lambda_0, \cdots, \lambda_n$ and the corresponding $\hat{r}(\lambda)$ values $\hat{r}_0, \cdots, \hat{r}_n$, if $n \geq 1$, $\lambda_{n+1}$ is obtained by using the secant formula on $\lambda_{n-1}$ and $\lambda_{n}$:
\begin{equation}
    \lambda_{n+1} = \lambda_n - \hat{r}_{n-1} \cdot
    {\lambda_{n-1} - \lambda_n \over \hat{r}_{n-1} - \hat{r}_n} \,.
\end{equation}
If $n=1$, we use a small value $\Delta \lambda$ to perturb $\lambda_{n-1}$:
\begin{equation}
    \lambda_{1} = \lambda_0  + \Delta \lambda \,.
\end{equation}
The updates terminate if a sufficient amount of training steps have been taken. Further details on the funnel inference algorithm may be found in the Supplementary Information.

\subsection*{Multi-stage Training} \label{methods:multistage}

If we only employ methods described up to this point, we are unable to reach maximum equation residuals much below $10^{-8}$. To discover higher-precision solutions, we employ multi-stage training\cite{Wang-Lai:multistage-nn}, which divides the training process into different stages, each with a dedicated neural network, c.f. hierarchical training~\cite{trask2022hierarchical, howard2023stacked}.
In this study, we apply two-stage training to the stable and first unstable solutions of the CCF and IPM equations. This approach lets us reach machine precision, as shown in the Supplementary Information. In the second stage, we train a new neural network that specifically learns the error $\mathcal{R}^\text{stage-1}_k$ between the first-stage networks and the exact solution. The final solution is obtained by composing the output of the two networks. This approach offers two advantages over single stage training. \\
\\
{\bf 1. Linearized governing equation:} Despite the nonlinearity of the original equation, the governing equation for the second stage training is effectively linear. We can formalize this mathematically by expressing the exact solution, $\hat{\Phi}_\text{exact}$, as the sum of the first-stage solution, $\hat{\Phi}_\text{stage-1}$, and a small error correction term: $\hat{\Phi}_\text{exact} = \hat{\Phi}_\text{stage-1} + \epsilon \hat{\Phi}_\text{stage-2}$. Assuming that the first-stage training converges properly, the magnitude of the error $\epsilon$ should be much less than 1. Then $\hat{\Phi}_\text{stage-2}$ becomes the normalized error profile that the second-stage network is trained to approximate. Substituting this ansatz into a governing PDE, we arrive at the governing equation for the error profile $\hat{\Phi}_\text{stage-2}$:
\begin{equation}\label{eq:burg_stage2}
	-\epsilon \mathcal{D}_k[\hat{\Phi}_\text{stage-1}] \hat{\Phi}_\text{stage-2} + O(\epsilon^2) = \mathcal{R}^\text{stage-1}_k \,,
\end{equation}
where $\mathcal{D}_k[\hat{\Phi}_\text{stage-1}]$ is the linear operator obtained by expanding the residual $\mathcal{R}_k$ around the approximate solution $\hat{\Phi}_\text{stage-1}$. By neglecting the $O(\epsilon^2)$ and higher-order terms, this equation becomes linear with respect to the function $\hat{\Phi}_\text{stage-2}$. The right-hand side of Equation~\eqref{eq:burg_stage2} is the equation residual from the first stage, which now acts as a known source term for the linear equation governing second stage training.\\
 \\
{\bf 2. Informed network structure:}
Given the near-linearity of the second stage governing equations, the properties of the first-stage residuals $\mathcal{R}^\text{stage-1}_k$ can inform the optimal architecture for the second-stage network. As shown in Wang \& Lai\cite{Wang-Lai:multistage-nn}, the characteristic features of the error profile $\hat{\Phi}_\text{stage-2}$ are governed by the source term $\mathcal{R}^\text{stage-1}_k$. For example, the equation residual typically exhibits higher frequency content than the solution $\hat{\Phi}_\text{stage-1}$. This implies that the error function $\mathcal{R}^\text{stage-1}_k$ is also a high-frequency function. Leveraging this fact, we can embed the appropriate inductive bias into the second stage network. Instead of using a standard MLP, we employ a Fourier feature network, which incorporates a Fourier mapping layer $[\cos({\bf B}x), \sin({\bf B}x)]$ placed after the linear input layer. The weights {\bf B} of this layer are sampled from a Gaussian distribution ${\bf B}\sim \mathcal{N}(0, \sigma)$, where $\sigma$ indicates the standard deviation of the distribution. We directly set the hyperparameter $\sigma$ using the dominant frequency $f_d^{(e)}$ of the first-stage equation residuals $\mathcal{R}^\text{stage-1}_k$, namely $\sigma = 2\pi f_d^{(e)}$. This provides a principled way of designing a network that is inherently suited to learn the high-frequency error function.

\subsection*{Validation of $\lambda$}

Once we discover a numerical solution $\hat{\Phi}^*$ and the associated self-similar scaling parameter $\lambda^*$, we also validate the admissibility of $\lambda^*$. Ideally, if we are able to evaluate the function $R(\lambda)$ which returns the lowest possible maximum residual achievable with a numerical function $\hat{\Phi}_\theta$ for a given $\lambda$, that is,
\begin{equation}
    R(\lambda) = \min_\theta \left[ \text{max d0 residual}_\text{grid} (\hat{\Phi}_\theta, \lambda) \right] \,,
\end{equation}
the admissibility of $\lambda$ can be validated at this level by requiring that $\lambda^*$ is a local minimum of $R(\lambda)$. We refer to $R(\lambda)$ as the \textit{smoothness error} of $\lambda$ in this section.

\begin{figure}[t]
	\centering
	\includegraphics[width=\linewidth]{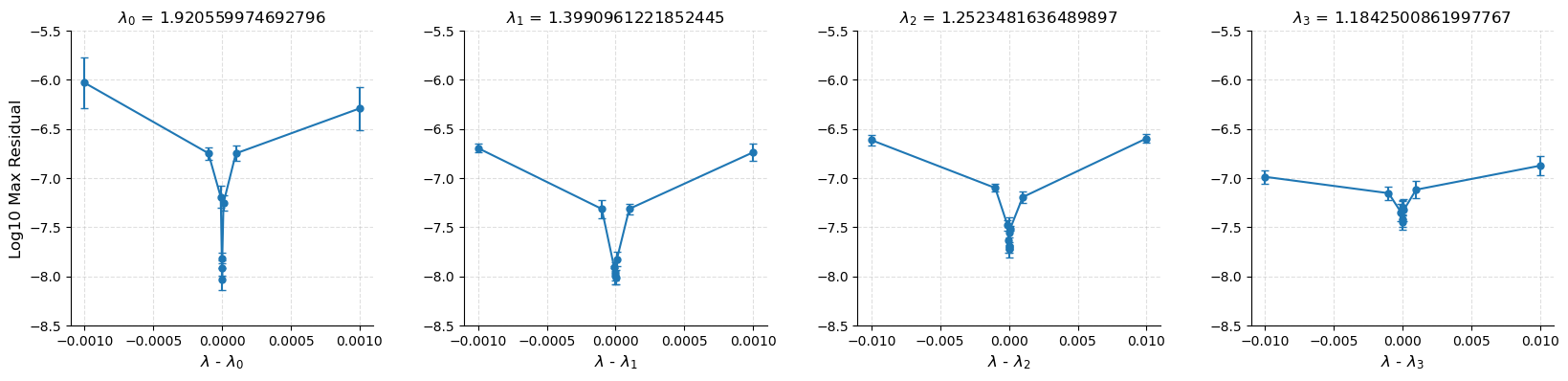}
	\caption{\small \textbf{Smoothness signal in $\lambda$.} Maximum residuals over the input space as a function of $\lambda - \lambda_i$ for admissible $\lambda_i$ of the Boussinesq equation with error bars. The plotted values and error bars are produced based on five experiments with varying seeds. Each subplot is centered around a different admissible lambda value, where $\lambda_0$ is the self-similar scaling parameter for the stable solution, and $\lambda_n$ with $n \geq 1$ is that for the $n$-th unstable solution. The x-axis shows the difference between the lambda value and the closest admissible lambda value ($\lambda - \lambda_i$), while the y-axis shows the mean and standard deviation of the log 10 maximum residual at 400k training steps. The plots show how the residuals increase as the lambda value deviates from $\lambda_i$ for each $i$, indicating that $\lambda_i$ is close to the local minimum of smoothness errors.}
	\label{fig:lamb_funnel_plot}
\end{figure}

Since we do not have the ability to evaluate $R$, we construct a proxy for the smoothness error by running the training pipeline that we have presented up to now, but with constant $\lambda$, and evaluating the maximum residual of the final solution obtained this way. In practice, this is done for a set number of neural networks initialized with different seeds. We take the geometric mean of the maximum residuals of the final solutions to produce the proxy function $\tilde{R}(\lambda)$ for $R(\lambda)$. We also compute the standard deviation of the metrics in log space to measure the uncertainty $\tilde{E}(\lambda)$ of this value. We construct a small exponentially spaced grid centered at $\lambda^*$, evaluate $\tilde{R}(\lambda)$ at those grid points and plot them against the $\lambda$ values. This is done for the $\lambda$ values discovered for the Boussinesq equation in figure \ref{fig:lamb_funnel_plot}. Note that the plots have a funnel shape centered around the discovered $\lambda$ values.

\begin{figure}[t]
	\centering
	\includegraphics[width=\linewidth]{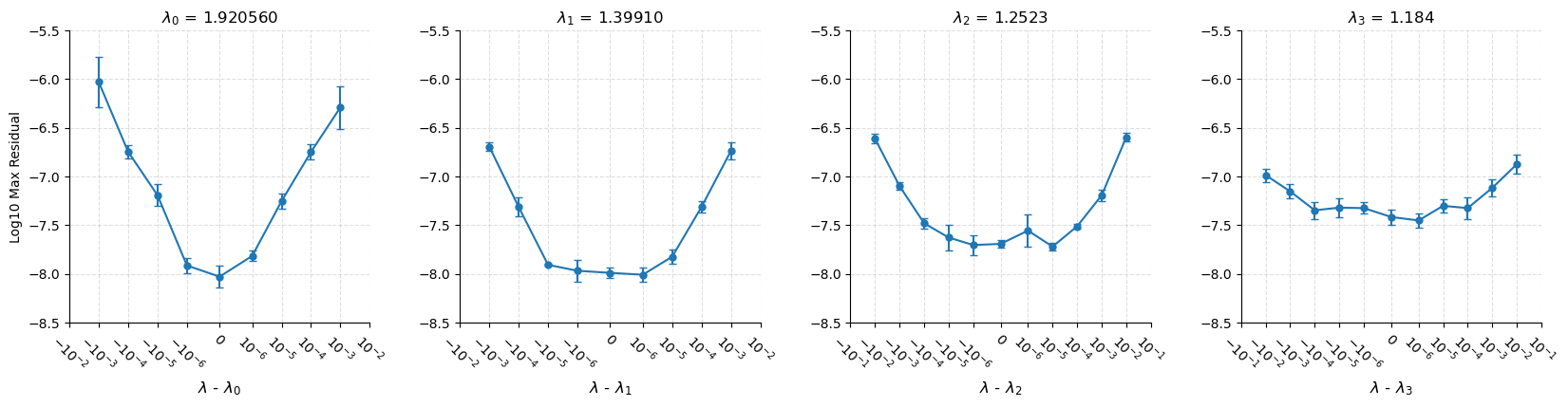}
	\caption{\small \textbf{Funnel plots for validation.} Figure \ref{fig:lamb_funnel_plot} with the x-axis in symlog-scale. The $\lambda_i$ values displayed now have been truncated to their significant digits based on the width of the basin of each plot.}
	\label{fig:lamb_funnel_plot_log_log}
\end{figure}

In figure \ref{fig:lamb_funnel_plot_log_log} we take the plots of figure \ref{fig:lamb_funnel_plot} and scale the $\lambda$-axis to be in symlog-scale. We find that the log-symlog plot of smoothness signal has a bowl shape due to the numerical nature of the evaluation method. In other words, there exists a limit to the resolution of $\lambda$ at which this method can discern admissibility.

Based on these observations, we use these funnel plots as a method to validate admissibility of a $\lambda^*$ value, and to determine its significant digits. First, we require that $\lambda^*$ lie at a basin of the bowl of the zoomed-in funnel plot in order to be considered valid. Furthermore, we interpret the width of the basin to be the resolution at which admissibility can be determined. We hence identify the last digit of the numerically attained admissible scaling parameter by the minimum perturbation that pushes the $\lambda$ value out of the basin of the bowl.

\subsection*{Spectrum of Linearization}

The self-similar profiles constructed are steady-state solutions to their respective governing equations when viewed in an appropriate rescaled or co-moving coordinate system. The linear stability of such a steady state, which we denote as $\Phi_\lambda$, is determined by the spectral properties of the operator $\mathcal{D}[\Phi_\lambda]$ derived by linearizing the system around $\Phi_\lambda$. This leads to an eigenvalue problem of the form:
$$\mathcal{D}[\Phi_\lambda] \Psi = \mu \Psi \,.$$
Here, the eigenfunction problem is restricted to $\Psi$ that lie within the same symmetry class as the original solution.

A profile $\Phi_\lambda$ is considered \textit{linearly stable} if the portion of its spectrum with a non-negative real part, $\Re(\mu) \ge 0$, consists solely of the trivial eigenvalues arising from the symmetries of the problem. If the spectrum contains any other eigenvalues in this region, the solution is considered \textit{linearly unstable}. The number of these non-trivial eigenvalues $\mu_i$ with a non-negative real part, $\Re(\mu_i) \geq 0$, is referred to as the \textit{order of instability}. For a computer-assisted proof to be feasible, it is desirable that the spectrum in the right-half $\mu$-plane consists of a finite number of eigenvalues.

In order to analyze the stability of the new solutions discovered in this work, we learn the non-negative real eigenvalues of the corresponding linearized operators, under the assumption that there exist eigenvalues with non-negative real part that lie on the real axis. We again employ PINNs to solve the eigenvalue problem, optimizing the eigenvalue $\mu$ simultaneously with the neural network representing the eigenfunction $\Psi$. The methods used to discover the solution of the self-similar equations themselves are largely transferred over to discovering the solutions to the eigenvalue problem. The PINNs are now trained to minimize eigenvalue equation residuals. The specific linearized operators, as well as details on the eigenfunction training setup for each problem can be found in the Supplementary Information.

\section*{Acknowledgements}

C-Y. Lai acknowledges the NSF  via grant  DMS-2245228. T. Buckmaster and Y. Wang were supported in part by the NSF grant DMS-2243205 and the Simons Foundation Mathematical and Physical Sciences collaborative grant `Wave Turbulence'. T. Buckmaster was also supported by the NSF grant DMS-2244879. Y. Wang's research was additionally supported by a gift to Buckmaster from Google as well as an NYU Postdoctoral Research and Professional Development Support Grant. Both T. Buckmaster and Y. Wang  gratefully acknowledge Google for providing computing credits and subscriptions, as well as   Stanford Research Computing and New York University High Performance Computing
for providing computational resources and support. G. Cao-Labora has been supported by NSF under grant DMS-2052651 during part of this work. G. Cao-Labora and J. G\'omez-Serrano have also been partially supported by the MICINN research grant number PID2021–125021NA–I00. J. G\'omez-Serrano was partially supported by NSF Grants DMS-2245017, DMS-2247537 and DMS-2434314; and by the Simons Foundation through a Simons Fellowship.

We thank A. P. Badia, J. Bausch,  C. Blundell, A. Botev, V. Dalibard, D. Hassabis, G. Holland, P.-S. Huang, Y. Li, A. Novikov, S. Nowozin, A. Pritzel, D. Rezende and P. Wirnsberger for their contributions, sponsorship, advice and helpful discussions. 

This work represents a collaboration between academic institutions and Google DeepMind.

\end{document}